\documentclass[12pt]{article}

\usepackage{amsmath,amsfonts,amssymb}

\def\R{{\bf R}}

\newtheorem{theorem}{Theorem}[section]
\newtheorem{lemma}[theorem]{Lemma}

\newtheorem{proposition}[theorem]{Proposition}

\newtheorem{remark}[theorem]{Remark}
\newtheorem{example}[theorem]{Example}

\begin{document}

\title{\vspace*{-1.5cm}
Approximations  of Stochastic Partial Differential Equations}

\author{ Giulia Di Nunno$^{1}$ \and
Tusheng Zhang$^{2}$}

\footnotetext[1]{\ Department of Mathematics, University of
Oslo, N-0316 Oslo, Norway}
\footnotetext[2]{\ School
of Mathematics, University of Manchester, Oxford Road, Manchester
M13 9PL, England, U.K. Email: tusheng.zhang@manchester.ac.uk}

\maketitle

\begin{abstract}
In this paper we show that solutions of stochastic partial differential equations driven by Brownian motion can be approximated by stochastic partial differential equations forced by pure jump noise/random kicks. Applications to stochastic Burgers equations are discussed.
\end{abstract}


\section{Introduction}

Stochastic evolution equations and stochastic partial differential
equations (SPDEs)  are of great interest to many
people. There exists a great amount of literature on the subject,
see, for example the monographs [DZ], [C]. \vskip 0.3cm
 \noindent In this paper, we consider the following stochastic evolution equation:
\begin{eqnarray}\label{eq:1.1}
dY_t&=&-AY_t dt+[b_1(Y_t)+b_2(Y_t)] dt +\sigma (Y_t) dB_t, \\
Y_0&=&h\in H,
\end{eqnarray}
in the framework of a Gelfand triple :
\begin{equation}
V\subset H\cong H^{\ast} \subset V^{\ast},
\end{equation}
where $H, V$ are Hilbert spaces, $A$ is the infinitesimal generator
of a strongly continuous semigroup, $b_1, \sigma$ are measurable
mappings from $H$ into $H$, $b_2$ is a measurable
mappings from $H$ into $V^*$(the dual of $V$), $B_t, t\geq 0$ is a Brownian motion.
The solutions are considered to be weak
solutions (in the PDE sense) in the space  $V$ and not as mild
solutions in $H$ as is more common in the literature. The stochastic
evolution equations of this type driven by  Wiener processes were
first studied  by in [P]
and subsequently in [KR].
For stochastic equations with general Hilbert space valued semimartingales
replacing the Brownian motion we refer to [GK1], [GK2], [G] and \cite{RZ}.

\vskip 0.3cm
 The aim of this paper is to study the approximations of stochastic evolution equations of the above type
 by solutions of stochastic evolution equations driven by pure jump processes, namely forced
 by random kicks. One of the motivations is to shine some light on numerical simulations of SPDEs driven by pure jump noise.  To include interesting applications, the drift of the equation (\ref{eq:1.1}) will consist of a ``good'' part $b_1$ and a ``bad'' part $b_2$. The crucial step of obtaining the approximation is to establish the tightness of the approximating equations in the space of Hilbert space-valued right continuous paths with left limits. This is tricky  because of the nature of the infinite dimensions and weak assumptions
 on the drift $b_2$. We first obtain the approximations assuming the diffusion coefficient $\sigma$ takes values
 in the smaller space $V$ and then remove the restriction by another layer of approximations. As far as we are aware of, this is the first paper to consider such approximations for SPDEs. The approximations of small jump L\`e{}vy processes were considered in \cite{AR}. Robustness of solutions of stochastic differential equations replacing small jump Levy processes by Brownian motion was discussed in \cite{BDK} and \cite{DKV}.
 \vskip 0.4cm
 The rest of the paper is organized as follows. In Section 2 we lay down the precise framework. The main part is Section 3 where approximations are established and applications to stochastic Burgers equations are discussed.

\section{Framework}

\setcounter{equation}{0}

Let $V$, $H$ be two separable Hilbert spaces such that $V$ is
continuously,  densely imbedded in $H$. Identifying $H$ with its dual
we have
\begin{equation}
V\subset H\cong H^{\ast} \subset V^{\ast},
\end{equation}
where $V^{\ast}$ stands for the topological dual of $V$. We assume that the imbedding $V\subset H$ is compact. Let $A$ be
a self-adjoint operator on the Hilbert space $H$ satisfying the
following coercivity hypothesis: There exist constants $\alpha_0 >0$, $\alpha_1 >0$
and $\lambda_0 \geq 0$ such that
\begin{equation}\label{2.1}
\alpha_0 ||u||_V^2\leq 2\langle Au,u\rangle +\lambda_0 |u|_H^2 \leq \alpha_1 ||u||_V^2\qquad
\hbox{for all $u\in V$}.
\end{equation}
$\langle Au,u\rangle =Au(u)$ denotes the action of $Au\in
V^{\ast}$ on $u\in V$.

\noindent We remark that $A$ is generally not bounded as an operator
from $H$ into $H$. Let $(\Omega, {\cal F}, P)$ be a probability
space equipped with a filtration $\{{\cal F}_t\}$ satisfying the
usual conditions. Let $\{B_t, t\geq 0\}$ be a real-valued ${\cal
F}_t$- Brownian motion, $\nu (dx)$  a $\sigma$-finite measure on the
measurable space $(R_0, {\cal B}(R_0))$, where $R_0=R\setminus \{0\}$. Let
$p=(p(t)), t\in D_p$ be a stationary  ${\cal F}_t$-Poisson point
process on $R_0$ with characteristic measure $\nu$. See [IW] for the
details on Poisson point processes. Denote by $ N(dt, dx)$ the
Poisson counting measure associated with $p$, i.e.,
$N(t,A)=\sum_{s\in D_p, s\leq t}I_A(p(s))$. Let
$\tilde{N}(dt,dx):=N(dt, dx)-dt\nu (dx)$ be the compensated Poisson
random measure. Let $b_1$, $\sigma$ be measurable mappings  from $H$
into  $ H$,  and  $b_2(\cdot)$  a measurable mapping from $H$
into $V^*$.
 Denote by $D([0,T], H)$ the space of all c\`a{}dl\`a{}g paths from
$[0,T]$ into $H$. Consider the  stochastic evolution  equations:
\begin{eqnarray}\label{2.2}
dX_t&=&-AX_t dt+[b_1(X_t)+b_2(X_t)] dt +\sigma (X_t) dB_t, \\
Y_0&=&h\in H.
\end{eqnarray}
\begin{eqnarray}\label{2.3}
X_t^{\varepsilon}&=&h-\int_0^tAX_s^{\varepsilon} ds+\int_0^t[b_1(X_s^{\varepsilon})+b_2(X_s^{\varepsilon})] ds \nonumber\\ &&+\frac{1}{\alpha(\epsilon)}\int_0^t\int_{|x|\leq \varepsilon}\sigma(X_{s-}^{\varepsilon})x\tilde{N}(ds,dx).
\end{eqnarray}

 Introduce the following conditions:
\vskip 0.2cm
\noindent {\bf (H.1)} There exists a constant $C<\infty$ such that
\begin{eqnarray}\label{2.5}
&& |b_1(y_1)-b_1(y_2)|_H^2 +|\sigma (y_1)-\sigma (y_2)|_H^2 \nonumber \\
&\leq &C |y_1-y_2|_H^2, \quad\quad \mbox{for all} \quad y_1, y_2\in H.
\end{eqnarray}
\vskip 0.2cm
\noindent {\bf (H.2)} $b_2(\cdot)$ is a mapping from $V$ into $V^*$ that satisfies

(i) $<b_2(u), u>=0$ for $u\in V$,

(ii) There exist constants $C_1$, $C_2$, $\beta<\frac{1}{2}$  such that
\begin{eqnarray}\label{2.6}
 &&<b_2(y_1)-b_2(y_2), y_1-y_2>\nonumber\\
&\leq &\beta \alpha_0 ||y_1-y_2||_V^2+C_1|y_1-y_2|_H^2(||y_1||_V^2+||y_2||_V^2)\nonumber\\
&& \quad \mbox{for all} \quad  y_1, y_2\in V,
\end{eqnarray}

(iii) There exists a constant $0<\gamma<1$ such that $||b_2(u)||_{V^*}\leq C_1 |u|_H^{2-\gamma}||u||_V^{\gamma}$ for $u\in V$.

\vskip 0.4cm
Under the assumptions (H.1) and (H.2), it is known that equations (\ref{2.2}), (\ref{2.3})
admit unique solutions.
\vskip 0.3cm

We finish this section with two
examples.
\begin{example}
\rm \ Let $D$ be a bounded domain in $\R^d$. Set $H=L^2(D)$. Let
$V=H_0^{1,2}(D)$ denote the Sobolev space of order one with homogenous boundary
conditions.
Denote by $a(x)=(a_{ij}(x))$  a symmetric matrix-valued function on $D$
satisfying  the uniform ellipticity condition:
$$
\frac{1}{c} I_d \leq a(x)\leq c I_d\qquad \hbox{for some constant
$\;c\in(0,\infty)$}.
$$
Define
$$
Au=-div(a(x)\nabla u(x)).
$$
Then (\ref{2.1}) is fulfilled for $(H,V,A)$.
\end{example}

\begin{example}
\rm \ Let $
Au=-\Delta_{\alpha},
$
where  $\Delta_{\alpha}$ denotes  the generator of a symmetric
$\alpha$-stable process in $R^d$, $0<\alpha \leq 2$.
$\Delta_{\alpha}$ is called the fractional Laplace operator. It is well known that the Dirichlet
form associated with $\Delta_{\alpha}$ is given by
$${\cal E}(u,v)=K(d,\alpha )\int \int_{R^d\times R^d} \frac{(u(x)-u(y))(v(x)-v(y))}{|x-y|^{d+\alpha}}\,dxdy,$$
$$D({\cal E})=\{u\in L^2(R^d):\quad \int \int_{R^d\times R^d} \frac{|u(x)-u(y)|^2}{|x-y|^{d+\alpha}}\,dxdy<\infty \},$$
where $K(d,\alpha )=\alpha 2^{\alpha
-3}\pi^{-\frac{d+2}{2}}sin(\frac{\alpha\pi}{2})\Gamma
(\frac{d+\alpha}{2})\Gamma (\frac{\alpha}{2})$. To study equation
(2.7), we choose  $H=L^2(\R^d)$, and $V=D({\cal E})$ with the inner
product $<u,v>={\cal E}(u,v)+(u,v)_{L^2(R^d)}$.
Then (\ref{2.1}) is fulfilled for $(H,V,A)$. See [FOT] for details about
the fractional Laplace operator.
\end{example}

\section{Approximations of SPDEs by pure jump type SPDEs}
Set
$$\alpha(\epsilon)=\left (\int_{\{|x|\leq \epsilon\}}x^2\nu(dx)\right )^{\frac{1}{2}}$$
Consider the following SPDE driven by pure jump noise:
\begin{eqnarray}\label{5.1}
X_t^{\varepsilon}&=&h-\int_0^tAX_s^{\varepsilon} ds+\int_0^t[b_1(X_s^{\varepsilon})+b_2(X_s^{\varepsilon})] ds \nonumber\\ &&+\frac{1}{\alpha(\epsilon)}\int_0^t\int_{|x|\leq \varepsilon}\sigma(X_{s-}^{\varepsilon})x\tilde{N}(ds,dx).
\end{eqnarray}
Let $X$ denote the solution to the SPDE:
\begin{eqnarray}\label{5.2}
X_t&=&h-\int_0^tAX_s ds+\int_0^t[b_1(X_s)+b_2(X_s)] ds +\int_0^t\sigma(X_{s})dB_s.
\end{eqnarray}
 Denote by $\mu_{\varepsilon}$, $\mu$ respectively the laws of $X^{\varepsilon}$ and $X$ on the spaces $D([0,T], H)$ and $C([0,T],H)$. Consider the following conditions:
\vskip 0.3cm
\noindent{\bf (H.3)} There exists a sequence of mappings $\sigma_n(\cdot): H\rightarrow V$ such that

(i) $|\sigma_n(y_1)-\sigma_n(y_2)|_H\leq c |y_1-y_2|_H$, where $c$ is a constant independent of $n$,

(ii) $|\sigma_n(y)-\sigma(y)|_H\rightarrow 0$ uniformly on bounded subsets of $H$.
\vskip 0.3cm
\begin{remark}
In most of the cases, one simply chooses $\sigma_n$ to be the finite dimensional projection of $\sigma$ into the space $V$.
\end{remark}
\vskip 0.3cm
\noindent{\bf (H.3)$^{\prime}$} The mapping $\sigma(\cdot)$ takes  the space $V$ into itself and  satisfies
 $||\sigma(y)||_V\leq c (1+||y||_V)$ for some constant $c$.
 \vskip 0.3cm
\noindent{\bf(H.4)} There exists an orthonormal basis $\{ e_k, k\geq 1\}$ of $H$ such that $Ae_k=\lambda_ke_k$ and $0\leq \lambda_1\leq \lambda_2\leq ...\leq \lambda_n\rightarrow \infty$
as $n\rightarrow \infty$.

\vskip 0.4cm
We first prepare some preliminary results needed for the proofs of the main theorems.
\vskip 0.3cm
The following estimate holds for $\{X^{\varepsilon},\varepsilon>0\}$.
\begin{lemma}
Let $X^{\varepsilon}$ be the solution of equation (\ref{5.1}). If $\frac{\varepsilon}{\alpha(\epsilon)}\leq C_0$ for some constant $C_0$, then we have for $p\geq 2$,
\begin{equation}\label{5.0}
\sup_{\varepsilon}\{E[\sup_{0\leq t\leq T}|X_t^{\varepsilon}|_H^p]+E[\left (\int_0^T||X_s^{\varepsilon}||_V^2ds\right )^{\frac{p}{2}}]\}<\infty.
\end{equation}
\end{lemma}
\vskip 0.2cm \noindent {\bf Proof}. We prove the lemma for $p=4$. Other cases are similar. In view of the assumption (H.2), by Ito's formula, we have
\begin{eqnarray}\label{5.01}
&&|X_t^{\varepsilon}|_H^2\nonumber\\
&=&|h|_H^2 -2\int_0^t<X_s^{\varepsilon}, AX_s^{\varepsilon}>ds+2 \int_0^t<X_s^{\varepsilon}, b_1(X_s^{\varepsilon})>ds\nonumber\\
&+&\int_0^t\int_{|x|\leq \varepsilon}\left ( |\frac{1}{\alpha(\epsilon)}\sigma(X_{s-}^{\varepsilon})x|_H^2+2<X_{s-}^{\varepsilon}, \frac{1}{\alpha(\epsilon)}\sigma(X_{s-}^{\varepsilon})x>\right ) \tilde{N}(ds,dx)\nonumber\\
&+&  \int_0^t\int_{|x|\leq \varepsilon}|\frac{1}{\alpha(\epsilon)}\sigma(X_{s-}^{\varepsilon})x|_H^2 ds \nu (dx).
\end{eqnarray}
Let
$$ M_t=\int_0^t\int_{|x|\leq \varepsilon}\left ( |\frac{1}{\alpha(\epsilon)}\sigma(X_{s-}^{\varepsilon})x|_H^2+2<X_{s-}^{\varepsilon}, \frac{1}{\alpha(\epsilon)}\sigma(X_{s-}^{\varepsilon})x>\right ) \tilde{N}(ds,dx).$$
By Burkh\"older's inequality, for $t\leq T$,
 \begin{eqnarray}\label{5.02}
&&E[\sup_{0\leq u\leq t}|M_u|_H^2]\leq C E[[M, M]_t]\nonumber\\
&=& CE[ \int_0^t\int_{|x|\leq \varepsilon}\left ( |\frac{1}{\alpha(\epsilon)}\sigma(X_{s-}^{\varepsilon})x|_H^2+2<X_{s-}^{\varepsilon}, \frac{1}{\alpha(\epsilon)}\sigma(X_{s-}^{\varepsilon})x>\right )^2 {N}(ds,dx)]  \nonumber\\
&=& CE[ \int_0^t\int_{|x|\leq \varepsilon}\left ( |\frac{1}{\alpha(\epsilon)}\sigma(X_{s-}^{\varepsilon})x|_H^2+2<X_{s-}^{\varepsilon}, \frac{1}{\alpha(\epsilon)}\sigma(X_{s-}^{\varepsilon})x>\right )^2 ds\nu(dx)]\nonumber\\
&\leq & CE[ \int_0^t(1+|X_{s}^{\varepsilon}|_H^4) ds],
\end{eqnarray}
where the linear growth condition on $\sigma$ and the fact $\frac{\varepsilon}{\alpha(\epsilon)}\leq C_0$ have been used.
Use first (\ref{2.1}) and then square both sides of the resulting inequality to obtain from (\ref{5.01}) that
 \begin{eqnarray}\label{5.01-1}
&&|X_t^{\varepsilon}|_H^4 + \left (\int_0^t||X_s^{\varepsilon}||_V^2ds\right )^2 \nonumber\\
&\leq &C_T|h|_H^4 + C_T \int_0^t(1+|X_{s}^{\varepsilon}|_H^4) ds+C_TM_t^2.
\end{eqnarray}
 Take superemum over the interval $[0,t]$ in (\ref{5.01-1}), use (\ref{5.02}) to get
\begin{eqnarray}\label{5.03}
&&E[\sup_{0\leq s\leq t}|X_s^{\varepsilon}|_H^4]+E[\left (\int_0^t||X_s^{\varepsilon}||_V^2ds\right )^2]\nonumber\\
&\leq & C|h|_H^4 + CE[ \int_0^t(1+|X_{s}^{\varepsilon}|_H^4) ds].
\end{eqnarray}
Applying Gronwall's inequality proves the lemma. $\blacksquare$
\vskip 0.3cm
\begin{proposition}
Assume (H.1), (H.2), (H.3)$^{\prime}$,(H.4) and $\frac{\varepsilon}{\alpha(\epsilon)}\leq C_0$. Then the family $\{X^{\varepsilon}, \varepsilon >0\}$ is tight on the space $D([0,T], H)$.
\end{proposition}
\vskip 0.2cm \noindent {\bf Proof}. Write
\begin{equation}\label{5.3}
Y_{t}^{\varepsilon}=\frac{1}{\alpha(\epsilon)}\int_0^t\int_{|x|\leq \varepsilon}\sigma(X_{s-}^{\varepsilon})x\tilde{N}(ds,dx),
\end{equation}
and set $Z_{t}^{\varepsilon}=X_{t}^{\varepsilon}-Y_{t}^{\varepsilon}$. It suffices to prove that both $\{Y^{\varepsilon}, \varepsilon >0\}$ and $\{Z^{\varepsilon}, \varepsilon >0\}$ are tight. This is done in two steps.
\vskip 0.3cm
\noindent{Step 1}. Prove that $\{Y^{\varepsilon}, \varepsilon >0\}$ is tight.
\vskip 0.3cm
In view of the assumptions on $\sigma$, we have $Y^{\varepsilon}\in D([0,T], V)$. Since the imbedding $V\subset H$ is compact, according to Theorem 3.1 in \cite{J}, it is sufficient to show that for every $e\in H$, $\{ <Y^{\varepsilon},e>, \varepsilon>0\}$ is tight in $D([0,T], R)$. Note that
\begin{eqnarray}\label{5.4}
&&\sup_{\varepsilon}E[\sup_{0\leq t\leq T}<Y_t^{\varepsilon},e>^2]\leq  \sup_{\varepsilon}E[\sup_{0\leq t\leq T}|Y_t^{\varepsilon}|_H^2]\nonumber\\
&\leq &C\sup_{\varepsilon}\frac{1}{\alpha(\epsilon)^2}E[\int_0^T\int_{|x|\leq \varepsilon}|\sigma(X_{s}^{\varepsilon})|_H^2x^2\nu(dx)ds]\nonumber\\
&=&C\sup_{\varepsilon}E[\int_0^T|\sigma(X_{s}^{\varepsilon})|_H^2 ds]<\infty,
\end{eqnarray}
and for any stoping times $\tau_{\varepsilon}\leq T$ and any positive constants $\delta_{\varepsilon}\rightarrow 0$
we have
\begin{eqnarray}\label{5.5}
&&E[|<Y_{\tau_{\varepsilon}} ^{\varepsilon},e>-<Y_{\tau_{\varepsilon}+\delta_{\varepsilon}} ^{\varepsilon},e>|^2]\leq \frac{1}{\alpha(\epsilon)^2}E[\int_{\tau_{\varepsilon}}^{\tau_{\varepsilon}+\delta_{\varepsilon}} \int_{|x|\leq \varepsilon}|\sigma(X_{s}^{\varepsilon})|_H^2x^2\nu(dx)ds]\nonumber\\
&\leq& C\delta_{\varepsilon}\sup_{\varepsilon} E[(1+\sup_{0\leq t\leq T}|X_{t}^{\varepsilon}|_H^2)]\rightarrow 0,
\end{eqnarray}
as $\varepsilon\rightarrow 0$.
By Theorem 3 in \cite{J}, (\ref{5.4}) and (\ref{5.5}) yields the tightness of $<Y^{\varepsilon},e>, \varepsilon>0$.
\vskip 0.3cm
\noindent{Step 2}. Prove that $\{Z^{\varepsilon}, \varepsilon >0\}$ is tight.
\vskip 0.3cm
It is easy to see that $Z^{\varepsilon}$ satisfies the equation:
\begin{eqnarray}\label{5.6}
Z_t^{\varepsilon}&=&h-\int_0^tAZ_s^{\varepsilon} ds-\int_0^tAY_s^{\varepsilon} ds\nonumber\\
&&+\int_0^tb_1(Z_s^{\varepsilon}+Y_s^{\varepsilon})ds+\int_0^tb_2(Z_s^{\varepsilon}+Y_s^{\varepsilon})ds.
\end{eqnarray}
Recall $\{e_k, k\geq 1\}$ is the othonormal basis of $H$ consisting of eigenvectors of $A$. We have
\begin{eqnarray}\label{5.7}
&&<Z_t^{\varepsilon},e_k>\nonumber\\
&=&<h,e_k>-\lambda_k\int_0^t<Z_s^{\varepsilon},e_k> ds-\lambda_k\int_0^t<Y_s^{\varepsilon},e_k> ds\nonumber\\
&& +\int_0^t<b_1(Z_s^{\varepsilon}+Y_s^{\varepsilon}), e_k> ds+\int_0^t<b_2(Z_s^{\varepsilon}+Y_s^{\varepsilon}), e_k> ds.
\end{eqnarray}
By Corollary 3 in \cite{J}, to obtain the tightness of $\{Z^{\varepsilon}, \varepsilon >0\}$ we need to show

(i).  $\{<Z^{\varepsilon}, e_k>,  \varepsilon >0\}$ is tight in $D([0,T], R)$ for every $k$,

(ii). for any $\delta>0$,
\begin{eqnarray}\label{5.8}
\lim_{N\rightarrow \infty }\sup_{\varepsilon}P(\sup_{0\leq t\leq T}R_N^{\varepsilon}(t)>\delta )=0,
\end{eqnarray}
where
$$R_N^{\varepsilon}(t)=\sum_{k=N}^{\infty}<Z^{\varepsilon}_t, e_k>^2.$$
The proof of (i) is similar to that of the tightness of $<Y^{\varepsilon},e>, \varepsilon>0$. It is omitted. Let us prove (ii). By the chain rule, it follows that
\begin{eqnarray}\label{5.9}
<Z_t^{\varepsilon},e_k>^2&=&<h,e_k>^2-2\lambda_k\int_0^t<Z_s^{\varepsilon},e_k>^2 ds-2\lambda_k\int_0^t<Y_s^{\varepsilon},e_k> <Z_s^{\varepsilon},e_k> ds\nonumber\\
&& +2\int_0^t<b_1(Z_s^{\varepsilon}+Y_s^{\varepsilon}), e_k> <Z_s^{\varepsilon},e_k> ds\nonumber\\
&&+2\int_0^t<b_2(Z_s^{\varepsilon}+Y_s^{\varepsilon}), e_k> <Z_s^{\varepsilon},e_k> ds.
\end{eqnarray}
By the variation of constants formula, we have
\begin{eqnarray}\label{5.10}
<Z_t^{\varepsilon},e_k>^2&=&e^{-2\lambda_kt}<h,e_k>^2-2\lambda_k\int_0^te^{-2\lambda_k(t-s)}<Y_s^{\varepsilon},e_k> <Z_s^{\varepsilon},e_k> ds\nonumber\\
&& +2\int_0^te^{-2\lambda_k(t-s)}<b_1(Z_s^{\varepsilon}+Y_s^{\varepsilon}), e_k> <Z_s^{\varepsilon},e_k> ds\nonumber\\
&&+2\int_0^te^{-2\lambda_k(t-s)}<b_2(Z_s^{\varepsilon}+Y_s^{\varepsilon}), e_k> <Z_s^{\varepsilon},e_k> ds.
\end{eqnarray}
Hence
\begin{eqnarray}\label{5.11}
&&R_N^{\varepsilon}(t)=\sum_{k=N}^{\infty}<Z_t^{\varepsilon},e_k>^2\nonumber\\
&=&\sum_{k=N}^{\infty}e^{-2\lambda_kt}<h,e_k>^2-2\int_0^t\sum_{k=N}^{\infty}\lambda_ke^{-2\lambda_k(t-s)}<Y_s^{\varepsilon},e_k> <Z_s^{\varepsilon},e_k> ds\nonumber\\
&& +2\int_0^t\sum_{k=N}^{\infty}e^{-2\lambda_k(t-s)}<b_1(Z_s^{\varepsilon}+Y_s^{\varepsilon}), e_k> <Z_s^{\varepsilon},e_k> ds\nonumber\\
&& +2\int_0^t\sum_{k=N}^{\infty}e^{-2\lambda_k(t-s)}<b_2(Z_s^{\varepsilon}+Y_s^{\varepsilon}), e_k> <Z_s^{\varepsilon},e_k> ds\nonumber\\
&=:& I_N^{(1)}(t)+I_N^{(2)}(t)+I_N^{(3)}(t)+I_N^{(4)}(t).
\end{eqnarray}
Obviously
\begin{equation}\label{5.12}
I_N^{(1)}(t)\leq \sum_{k=N}^{\infty}<h,e_k>^2\rightarrow 0,
\end{equation}
as $N\rightarrow \infty$.
For the third term on the right side of (\ref{5.11}), we have
\begin{eqnarray}\label{5.13}
|I_N^{(3)}(t)|&\leq & 2\int_0^te^{-2\lambda_N(t-s)}\sum_{k=N}^{\infty}|<b_1(Z_s^{\varepsilon}+Y_s^{\varepsilon}), e_k> <Z_s^{\varepsilon},e_k>| ds\nonumber\\
&\leq &2\int_0^te^{-2\lambda_N(t-s)}ds (\sup_{0\leq s\leq T}|Z_s^{\varepsilon}|_H)(\sup_{0\leq s\leq T}|b_1(Z_s^{\varepsilon}+Y_s^{\varepsilon})|_H)\nonumber\\
&\leq& C\frac{1}{\lambda_N} \left (1+\sup_{0\leq s\leq T}|Z_s^{\varepsilon}|_H^2+\sup_{0\leq s\leq T}|Y_s^{\varepsilon}|_H^2\right ).
\end{eqnarray}
Hence,
\begin{eqnarray}\label{5.14}
&& \sup_{\varepsilon}E[\sup_{0\leq t\leq T}|I_N^{(3)}(t)|]\nonumber\\
&\leq & C\frac{1}{\lambda_N} \left (1+\sup_{\varepsilon}E[\sup_{0\leq s\leq T}|Z_s^{\varepsilon}|_H^2]+\sup_{\varepsilon}E[\sup_{0\leq s\leq T}|Y_s^{\varepsilon}|_H^2]\right )\nonumber\\
&&\rightarrow 0, \quad \mbox{as}\quad N\rightarrow \infty.
\end{eqnarray}
Let us turn to $I_N^{(2)}(t)$. By H\"older's inequality,
\begin{eqnarray}\label{5.15}
|I_N^{(2)}(t)|&\leq & 2\int_0^t(\sum_{k=N}^{\infty}e^{-4\lambda_k(t-s)}\lambda_k<Z_s^{\varepsilon}, e_k>^2)^{\frac{1}{2}}\nonumber\\
&&\quad\quad \times  (\sum_{k=N}^{\infty}\lambda_k<Y_s^{\varepsilon}, e_k>^2)^{\frac{1}{2}}ds\nonumber\\
&\leq & 2\int_0^t(\sum_{k=N}^{\infty}e^{-4\lambda_k(t-s)}\lambda_k<Z_s^{\varepsilon}, e_k>^2)^{\frac{1}{2}}(<AY_s^{\varepsilon}, Y_s^{\varepsilon}>^2)^{\frac{1}{2}}ds\nonumber\\
&\leq &C(\sup_{0\leq s\leq T}||Y_s^{\varepsilon}||_V)\int_0^te^{-\lambda_N(t-s)}(\sum_{k=N}^{\infty}e^{-2\lambda_k(t-s)}\lambda_k<Z_s^{\varepsilon}, e_k>^2)^{\frac{1}{2}}ds\nonumber\\
&\leq &C(\sup_{0\leq s\leq T}||Y_s^{\varepsilon}||_V)(\sup_{0\leq s\leq T}|Z_s^{\varepsilon}|_H)\int_0^te^{-\lambda_N(t-s)}\frac{1}{\sqrt{t-s}}ds\nonumber\\
&\leq &C(\frac{1}{\sqrt{\lambda_N}}\int_0^{\infty}e^{-u}\frac{1}{\sqrt{u}}du)(\sup_{0\leq s\leq T}||Y_s^{\varepsilon}||_V)(\sup_{0\leq s\leq T}|Z_s^{\varepsilon}|_H).
\end{eqnarray}
 In view of the assumption (H.2), the last term on the right side of (\ref{5.11}) can be estimated as follows:
\begin{eqnarray}\label{3.19}
|I_N^{(4)}(t)|&=&|\int_0^t\sum_{k=N}^{\infty}e^{-2\lambda_k(t-s)}<b_2(X_s^{\varepsilon}),e_k> <Z^{\varepsilon}_s,e_k> ds|\nonumber\\
&=&|\int_0^t\sum_{k=N}^{\infty}e^{-2\lambda_k(t-s)} \sqrt{\lambda_0+\lambda_k}<(A+\lambda_0I)^{-\frac{1}{2}}b_2(X_s^{\varepsilon}),e_k> <Z_s^{\varepsilon},e_k> ds|\nonumber\\
&\leq& C \int_0^t\left (\sum_{k=N}^{\infty}e^{-4\lambda_k(t-s)} <(A+\lambda_0I)^{-\frac{1}{2}}b_2(X_s^{\varepsilon}),e_k>^2  \right )^{\frac{1}{2}}\nonumber\\
&&\quad\times \left (\sum_{k=N}^{\infty}(\lambda_0+\lambda_k) <Z_s^{\varepsilon},e_k>^2\right )^{\frac{1}{2}} ds\nonumber\\
&\leq& C \int_0^t||Z^{\varepsilon}_s||_V e^{-2\lambda_N(t-s)}\left (\sum_{k=N}^{\infty}<(A+\lambda_0I)^{-\frac{1}{2}}b_2(X_s^{\varepsilon}),e_k>^2  \right )^{\frac{1}{2}}ds \nonumber\\
&\leq& C \int_0^t||Z^{\varepsilon}_s||_V e^{-2\lambda_N(t-s)}||b_2(X_s^{\varepsilon})||_{V^*} ds\nonumber\\
&\leq& C \int_0^t||Z^{\varepsilon}_s||_V e^{-2\lambda_N(t-s)} |X_s^{\varepsilon}|_H^{2-\gamma} ||X_s^{\varepsilon}||_V^{\gamma}ds.
\end{eqnarray}
This yields that
\begin{eqnarray}\label{3.07}
|I_N^{(4)}(t)|
&\leq&C \sup_{0\leq s\leq T}|X_s^{\varepsilon}|_{H}^{2-\gamma}\int_0^t||Z_s^{\varepsilon}||_V e^{-2\lambda_N(t-s)}||X_s^{\varepsilon}||_{V}^{\gamma}ds\nonumber\\
&\leq&C \sup_{0\leq s\leq T}|X_s^{\varepsilon}|_{H}^{2-\gamma}\int_0^t e^{-2\lambda_N(t-s)}(||X_s^{\varepsilon}||_{V}^{1+\gamma}+ ||X_s^{\varepsilon}||_{V}^{\gamma}||Y_s^{\varepsilon}||_{V}) ds\nonumber\\
&\leq&C \sup_{0\leq s\leq T}|X_s^{\varepsilon}|_{H}^{2-\gamma}\left (\int_0^t e^{-\frac{4}{1-\gamma} \lambda_N(t-s)}ds\right )^{\frac{1-\gamma}{2}} \nonumber\\
&&\quad\quad \times \left (\int_0^T \left (||X_s^{\varepsilon}||_{V}^{1+\gamma}+ ||X_s^{\varepsilon}||_{V}^{\gamma}||Y_s^{\varepsilon}||_{V}\right )^{\frac{2}{1+\gamma}}ds\right )^{\frac{1+\gamma}{2}} \nonumber\\
&\leq&C (\frac{1}{\lambda_N})^{\frac{1-\gamma}{2}} \sup_{0\leq s\leq T}|X_s^{\varepsilon}|_{H}^{2-\gamma} \left (\int_0^T \left (||X_s^{\varepsilon}||_{V}^{1+\gamma}+ ||X_s^{\varepsilon}||_{V}^{\gamma}||Y_s^{\varepsilon}||_{V}\right )^{\frac{2}{1+\gamma}}ds\right )^{\frac{1+\gamma}{2}}\nonumber\\
&&
\end{eqnarray}
Hence,
\begin{eqnarray}\label{5.15-1}
&& \sup_{\varepsilon}E[\sup_{0\leq t\leq T}|I_N^{(4)}(t)|]\nonumber\\
&\leq &  C (\frac{1}{\lambda_N})^{\frac{1-\gamma}{2}}\sup_{\varepsilon}E\left [ \sup_{0\leq s\leq T}|X_s^{\varepsilon}|_{H}^{2-\gamma} \right.\nonumber\\
&&\left. \times\left (\int_0^T \left (||X_s^{\varepsilon}||_{V}^{1+\gamma}+ ||X_s^{\varepsilon}||_{V}^{\gamma}||Y_s^{\varepsilon}||_{V}\right )^{\frac{2}{1+\gamma}}ds\right )^{\frac{1+\gamma}{2}}\right ] \nonumber\\
&\leq &  C (\frac{1}{\lambda_N})^{\frac{1-\gamma}{2}}\sup_{\varepsilon}E\left [\sup_{0\leq s\leq T}|X_s^{\varepsilon}|_{H}^{2-\gamma}\right.\nonumber\\
 &&\left.\quad \times \left (\int_0^T \left (C||X_s^{\varepsilon}||_{V}^{2}+ c||Y_s^{\varepsilon}||_{V}^2\right )ds\right )^{\frac{1+\gamma}{2}}\right ] \nonumber\\
&&\rightarrow 0, \quad \mbox{as}\quad N\rightarrow \infty,
\end{eqnarray}
where we used the fact that
$$|ab|\leq C(|a|^p+|b|^q), \quad \frac{1}{p}+\frac{1}{q}=1.$$
Putting together (\ref{5.11})---(\ref{5.15-1}) and applying the Chebychev inequality  we obtain (\ref{5.8}).$\blacksquare$
\vskip 0.3cm
Let ${\cal D}$ denote the class of functions $ f\in C_b^3(H)$ that satisfy (i) $f^{\prime}(z)\in D(A)$ and $|Af^{\prime}(z)|_H\leq C(1+|z|_H)$
for some constant $C$, where $f^{\prime}(z)$ stands for the Frechet derivative of $f$, (ii) $f^{\prime\prime}$, $f^{\prime\prime\prime}$ are bounded.
\vskip 0.3cm

For $f\in {\cal D}$, define
\begin{eqnarray}\label{5.16}
&&L^{\varepsilon}f(z)=-<Af^{\prime}(z),z>+<b_1(z),f^{\prime}(z)>+<b_2(z),f^{\prime}(z)>\nonumber\\
&&+\int_{|x|\leq \varepsilon}\left [ f(z+\frac{1}{\alpha(\epsilon)}\sigma(z)x)-f(z)-<f^{\prime}(z),\frac{1}{\alpha(\epsilon)}\sigma(z)x>\right ]\nu(dx),\nonumber\\
&&
\end{eqnarray}
and
\begin{eqnarray}\label{5.17}
Lf(z)&=& -<Af^{\prime}(z),z>+<b_1(z),f^{\prime}(z)>+<b_2(z),f^{\prime}(z)>\nonumber\\
&&+\frac{1}{2}<f^{\prime\prime}(z)\sigma(z), \sigma(z)>.
\end{eqnarray}
\begin{lemma}
Assume $\lim_{\varepsilon\rightarrow 0}\frac{\varepsilon}{\alpha(\varepsilon)}=0$. For $f\in {\cal D}$, it holds that
\begin{eqnarray}\label{5.18}
L^{\varepsilon}f(z)\rightarrow Lf(z) \quad \mbox{uniformly on bounded subsets of }\quad H
\end{eqnarray}
as $\varepsilon\rightarrow 0$.
\end{lemma}
\vskip 0.2cm \noindent {\bf Proof}. Note that
$$f(y+w)-f(y)-<f^{\prime}(y),w>=\int_0^1d\alpha \int_0^{\alpha} <f^{\prime\prime}(y+\beta w)w,w>d\beta.$$
Thus
\begin{eqnarray}\label{5.19}
&& L^{\varepsilon}f(z)-Lf(z)\nonumber\\
&=&\int_{|x|\leq \varepsilon}\int_0^1d\alpha \int_0^{\alpha}d\beta <f^{\prime\prime}(z+\beta\frac{1}{\alpha(\epsilon)}\sigma(z)x)\frac{1}{\alpha(\epsilon)}\sigma(z)x, \frac{1}{\alpha(\epsilon)}\sigma(z)x >\nu(dx)\nonumber\\
&&\quad\quad\quad -\int_0^1d\alpha \int_0^{\alpha}d\beta <f^{\prime\prime}(z)\sigma(z), \sigma(z)>\nonumber\\
&=&\frac{1}{\alpha(\epsilon)^2}\int_{|x|\leq \varepsilon}x^2\nu(dx)\int_0^1d\alpha \int_0^{\alpha}d\beta \left [ <f^{\prime\prime}(z+\beta\frac{1}{\alpha(\epsilon)}\sigma(z)x)\sigma(z), \sigma(z) >\right. \nonumber\\
&&\quad\quad \left. -<f^{\prime\prime}(z)\sigma(z), \sigma(z)>\right ]
\end{eqnarray}
Hence, for $z\in B_N=\{ z\in H; |z|_H\leq N\}$ we have
\begin{eqnarray}\label{5.20}
&& |L^{\varepsilon}f(z)-Lf(z)|\nonumber\\
&\leq &C\frac{1}{\alpha(\epsilon)^2}\int_{|x|\leq \varepsilon}x^2\nu(dx)\int_0^1d\alpha \int_0^{\alpha}d\beta \beta\frac{1}{\alpha(\epsilon)}|\sigma(z)|_H|x| |\sigma(z)|_H^2\nonumber\\
&\leq & C_N \frac{\varepsilon}{\alpha(\epsilon)}\rightarrow 0,
\end{eqnarray}
uniformly on $B_N$ as $\varepsilon\rightarrow 0$, where we have used the local Lipschtiz continuity of $f^{\prime\prime}(z)$.
$\blacksquare$

\vskip 0.4cm
\begin{theorem}
Suppose (H.1), (H.2),(H.3)$^{\prime}$, (H.4) hold and  $\lim_{\varepsilon\rightarrow 0}\frac{\varepsilon}{\alpha(\varepsilon)}=0$. Then,
for any $T>0$, $\mu_{\varepsilon}$ converges weakly to $\mu$ on the space $D([0,T], H)$ equipped with the Skorokhod topology.
\end{theorem}
\vskip 0.3cm
\noindent {\bf Proof}. Since the mapping $\sigma$ takes values in the space $V$,
by Proposition 3.3, the family $\{ \mu_{\varepsilon}, \varepsilon>0\}$ is tight.  Let $ \mu_0$ be the weak limit of any convergent sequence $\{\mu_{\varepsilon_n}\}$ on the canonical space $(\Omega=D([0,T], H), {\cal F})$ as $\varepsilon_n\rightarrow 0$. We will show that $\mu_0=\mu$. Denote by $X_t(\omega)=w(t), \omega \in \Omega$ the coordinate process.  Set $J(X)=\sup_{0\leq s\leq T}|X_{s}-X_{s-}|_H$. Since
\begin{eqnarray}\label{5.21}
&&  E^{\mu_{\varepsilon}}[J(X)]=E[J(X^{\varepsilon})]\nonumber\\
&\leq &\frac{\varepsilon}{\alpha(\epsilon)}E[\sup_{0\leq s\leq T}|\sigma(X^{\varepsilon}_s)|_H]\nonumber\\
&\leq & C \frac{\varepsilon}{\alpha(\epsilon)}(1+ E[\sup_{0\leq s\leq T}|X^{\varepsilon}_s|_H])\rightarrow 0,
\end{eqnarray}
as $\varepsilon \rightarrow 0$, it follows from Theorem 13.4 in \cite{B} that $\mu_0$ is supported on the $C([0,T],H)$, the space of $H$-valued continuous functions on $[0,T]$. As a consequence, the finite dimensional distributions of $\mu_{\varepsilon_n}$ converge to that of $\mu_0$.
\vskip 0.3cm
Let $f\in {\cal D}$. By Ito's formula,
\begin{eqnarray}\label{5.22}
&& f(X_{t}^{\varepsilon})-f(h)-\int_0^tL^{\varepsilon}f(X_{s}^{\varepsilon})ds \nonumber\\
&=& \int_0^t\int_{|x|\leq \varepsilon}\{f(X_{s-}^{\varepsilon}+ \frac{1}{\alpha(\epsilon)}\sigma(X_{s-}^{\varepsilon})x)- f(X_{s-}^{\varepsilon})\}  \tilde{N}(ds,dx)
\end{eqnarray}
is a martingale. Hence, for any $s_0< s_1<...<s_n\leq s<t$ and $f_0, f_1, ...f_n\in C_b(H)$ it holds that
\begin{eqnarray}\label{5.23}
&& E^{\mu_{\varepsilon}}[\left(f(X_t)-f(X_s)-\int_s^tL^{\varepsilon}f(X_{u})du\right)f(X_{s_0})...f(X_{s_n})] \nonumber\\
&=& 0.
\end{eqnarray}
For any positive constant $M>0$, by Lemma 3.4 we have
\begin{equation}\label{5.24}
\lim_{n\rightarrow \infty}E^{\mu_{\varepsilon_n}}[\int_s^t|L^{\varepsilon_n}f(X_{u})-Lf(X_{u})|du, \sup_{0\leq u\leq T}|X_u|_H\leq M ]=0.
\end{equation}
On the other hand, in view of the assumptions on $f$ we have
  \begin{eqnarray}\label{5.25}
&&\sup_{n}E^{\mu_{\varepsilon_n}}[\int_s^t|L^{\varepsilon_n}f(X_{u})-Lf(X_{u})|du, \sup_{0\leq u\leq T}|X_u|_H> M ]\nonumber\\
&\leq& C \frac{1}{M}\sup_{n}E^{\mu_{\varepsilon_n}}[\sup_{0\leq u\leq T}|X_u|_H^3]\leq C^{\prime} \frac{1}{M}
\end{eqnarray}
Combining (\ref{5.24}) with (\ref{5.25}) we arrive at
\begin{equation}\label{5.26}
\lim_{n\rightarrow \infty}E^{\mu_{\varepsilon_n}}[\int_s^t|L^{\varepsilon_n}f(X_{u})-Lf(X_{u})|du]=0.
\end{equation}
By the weak convergence of $\mu_{\varepsilon_n}$ and  the convergence of finite distributions, it follows from (\ref{5.23}) and (\ref{5.26}) that
\begin{eqnarray}\label{5.27}
&& E^{\mu_{0}}[(f(X_t)-f(X_s)-\int_s^tLf(X_{u})du)f(X_{s_0})...f(X_{s_n})] \nonumber\\
&=&\lim_{n\rightarrow \infty}E^{\mu_{\varepsilon_n}}[(f(X_t)-f(X_s)-\int_s^tLf(X_{u})du)f(X_{s_0})...f(X_{s_n})] \nonumber\\
&=&\lim_{n\rightarrow \infty}E^{\mu_{\varepsilon_n}}[(f(X_t)-f(X_s)-\int_s^tL^{\varepsilon_n}
f(X_{u})du)f(X_{s_0})...f(X_{s_n})]\nonumber\\
&= &0.
\end{eqnarray}
Since $s_0< s_1<...<s_n\leq s<t$ are arbitrary, (\ref{5.27}) implies that for any $f\in {\cal D}$,
$$M_t^f=f(X_{t})-f(h)-\int_0^tLf(X_{s})ds, \quad  t\geq 0,$$
is a martingale under $\mu_0$.
In particular, let $f(z)=<e_k, z>$ and $f(z)=<e_k,z><e_j,z>$ respectively to obtain that under $\mu_0$
\begin{eqnarray}
M^k_t&:=&<e_k,X_t>-<e_k,h>+\int_0^t<Ae_k, X_s>ds-\int_0^t<b_1(X_s),e_k>ds\nonumber\\
&&\quad\quad -\int_0^t<b_2(X_s),e_k>ds
\end{eqnarray}
and
\begin{eqnarray}
&&M^{k,j}_t:=<e_k,X_t><e_j,X_t>-<e_k,h><e_j,h>\nonumber\\
&&+\int_0^t\{<Ae_k, X_s><e_j,X_s>+<Ae_j, X_s><e_k,X_s>\}ds\nonumber\\
&&-\int_0^t<b_1(X_s),e_k<e_j,X_s>+e_j<e_k,X_s> >ds\nonumber\\
&&-\int_0^t<b_2(X_s),e_k<e_j,X_s>+e_j<e_k,X_s> >ds\nonumber\\
&&-\int_0^t<\sigma(X_s),e_k><\sigma(X_s),e_j> ds
\end{eqnarray}
are martingales. This together with Ito's formula yields that
\begin{equation}\label{5.28}
<M^k, M^j>_t=\int_0^t<\sigma(X_s),e_k><\sigma(X_s),e_j> ds,
\end{equation}
where $<M^k, M^j>$ stands for the sharp bracket of  the two martingales.
Now by  Theorem 18.12 in \cite{K} (or  Theorem $7.1^{\prime}$ in \cite{IW}), there exists a probability space
$(\Omega^{\prime}, {\cal F}^{\prime}, P^{\prime})$ with a filtration ${\cal F}_t^{\prime}$ such that on the standard extension
$$(\Omega\times\Omega^{\prime},  {\cal F}\times {\cal F}^{\prime}, {\cal F}_t\times {\cal F}^{\prime}_t, \mu_0\times P^{\prime} )$$
of  $(\Omega, {\cal F}, {\cal F}_t, P)$ there exists a Brownian motion $B_t,t\geq 0$ such that
\begin{equation}\label{5.29}
M^k_t = \int_0^t<\sigma(X_s),e_k>dB_s,
\end{equation}
namely,
\begin{eqnarray}\label{5.29-1}
&&<e_k,X_t>-<e_k,h>\nonumber\\
&=&-\int_0^t<Ae_k, X_s>ds+\int_0^t<b_1(X_s),e_k>ds+\int_0^t<b_2(X_s),e_k>ds\nonumber\\
&&+\int_0^t<\sigma(X_s),e_k>dB_s
\end{eqnarray}
for any $k\geq 1$. Thus, under $\mu_0$, $X_t, t\geq 0$ is a weak solution (both in the probabilistic and in PDE sense)  of the SPDE:
$$X_t=h-\int_0^tAX_sds+\int_0^tb_1(X_s)ds+\int_0^tb_2(X_s)ds+\int_0^t\sigma(X_s)dB_s $$
By the uniqueness of the above equation, we conclude that $\mu_0=\mu$ completing the proof of the theorem.$\blacksquare$

\vskip 0.3cm
\begin{theorem}
Suppose (H.1), (H.2), (H.3) and (H.4) hold and  $\lim_{\varepsilon\rightarrow 0}\frac{\varepsilon}{\alpha(\varepsilon)}=0$. Then,
for any $T>0$, $\mu_{\varepsilon}$ converges weakly to $\mu$ on the space $D([0,T], H)$ equipped with the Skorokhod topology.
\end{theorem}
\vskip 0.3cm
\noindent {\bf Proof}. Let $ \sigma_n(\cdot)$ be the mapping specified in (H.3). Let $X^{n,\varepsilon}, X^n$ be the solutions of the SPDEs:
\begin{eqnarray}\label{5.30}
X_t^{n,\varepsilon}& = &h-\int_0^tAX_s^{n,\varepsilon} ds+\int_0^tb_1(X_s^{n,\varepsilon}) ds+\int_0^tb_2(X_s^{n,\varepsilon}) ds \nonumber\\ &&+\frac{1}{\alpha(\epsilon)}\int_0^t\int_{|x|\leq \varepsilon}\sigma_n(X_{s-}^{n,\varepsilon})x\tilde{N}(ds,dx).
\end{eqnarray}
\begin{eqnarray}\label{5.31}
X_t^n&=&h-\int_0^tAX_s^n ds+\int_0^tb_1(X_s^n) ds+\int_0^tb_2(X_s^n) ds\nonumber\\
&& +\int_0^t\sigma_n(X_{s}^n)dB_s.
\end{eqnarray}
We claim  that for any $\delta>0$,
\begin{eqnarray}\label{5.32}
\lim_{n\rightarrow \infty}\sup_{\varepsilon}P(\sup_{0\leq t\leq T}|X_t^{n,\varepsilon}-X_t^{\varepsilon}|>\delta  )=0.
\end{eqnarray}
\begin{eqnarray}\label{5.33}
\lim_{n\rightarrow \infty}P(\sup_{0\leq t\leq T}|X_t^{n}-X_t|^2>\delta )=0.
\end{eqnarray}
Let us only prove (\ref{5.32}). The proof of (\ref{5.33}) is simpler.  As the proof of (\ref{5.0}), we can show that
\begin{equation}\label{5.33-1}
\sup_n\sup_{\varepsilon}\{E[\sup_{0\leq t\leq T}|X_t^{n,\varepsilon}|_H^2]+E[\int_0^T||X_s^{n,\varepsilon}||_V^2ds]\}<\infty.
\end{equation}
\begin{equation}\label{5.33-2}
\sup_{n}\{E[\sup_{0\leq t\leq T}|X_t^{n}|_H^2]+E[\int_0^T||X_s^{n}||_V^2ds]\}<\infty.
\end{equation}
By Ito's formula, we have
\begin{eqnarray}\label{5.34}
&&e^{-\gamma \int_0^t(||X_s^{n,\varepsilon}||_V^2+||X_s^{\varepsilon}||_V^2)ds}|X_t^{n,\varepsilon}-X_t^{\varepsilon}|_H^2\nonumber\\
&=&-\gamma \int_0^t e^{-\gamma \int_0^s(||X_u^{n,\varepsilon}||_V^2+||X_u^{\varepsilon}||_V^2)du}
|X_s^{n,\varepsilon}-X_s^{\varepsilon}|_H^2(||X_s^{n,\varepsilon}||_V^2+||X_s^{\varepsilon}||_V^2)ds\nonumber\\
&& -2\int_0^te^{-\gamma \int_0^s(||X_u^{n,\varepsilon}||_V^2+||X_u^{\varepsilon}||_V^2)du}<X_s^{n,\varepsilon}-X_s^{\varepsilon}, A(X_s^{n,\varepsilon}-X_s^{\varepsilon})>ds\nonumber\\
&&+2 \int_0^te^{-\gamma \int_0^s(||X_u^{n,\varepsilon}||_V^2+||X_u^{\varepsilon}||_V^2)du}<X_s^{n,\varepsilon}-X_s^{\varepsilon}, b_1(X_s^{n,\varepsilon})-b_1(X_s^{\varepsilon})>ds\nonumber\\
&&+2 \int_0^te^{-\gamma \int_0^s(||X_u^{n,\varepsilon}||_V^2+||X_u^{\varepsilon}||_V^2)du}<X_s^{n,\varepsilon}-X_s^{\varepsilon}, b_2(X_s^{n,\varepsilon})-b_2(X_s^{\varepsilon})>ds\nonumber\\
&+&\int_0^t\int_{|x|\leq \varepsilon}e^{-\gamma \int_0^s(||X_u^{n,\varepsilon}||_V^2+||X_u^{\varepsilon}||_V^2)du}\left ( |\frac{1}{\alpha(\epsilon)}(\sigma_n(X_{s-}^{n,\varepsilon})x-\sigma(X_{s-}^{\varepsilon})x)|_H^2\right.\nonumber\\
&&\quad\quad\left.  +2<(X_s^{n,\varepsilon}-X_s^{\varepsilon}), \frac{1}{\alpha(\epsilon)}(\sigma_n(X_{s-}^{n,\varepsilon})x-\sigma(X_{s-}^{\varepsilon})x)>\right ) \tilde{N}(ds,dx)\nonumber\\
&+&  \int_0^t\int_{|x|\leq \varepsilon}e^{-\gamma \int_0^s(||X_u^{n,\varepsilon}||_V^2+||X_u^{\varepsilon}||_V^2)du}|\frac{1}{\alpha(\epsilon)}(\sigma_n(X_{s-}^{n,\varepsilon})x-\sigma(X_{s-}^{\varepsilon})x)|_H^2 ds \nu (dx)\nonumber\\
&:=& \sum_{k=1}^6 I_k^{n, \varepsilon}(t).
\end{eqnarray}
In view of the assumption (\ref{2.6}), we see that
\begin{eqnarray}\label{5.34-1}
&&I_1^{n, \varepsilon}(t)+I_2^{n, \varepsilon}(t)+I_4^{n, \varepsilon}(t) \nonumber\\
&\leq &-(1-2\beta )\alpha_0 \int_0^te^{-\gamma \int_0^s(||X_u^{n,\varepsilon}||_V^2+||X_u^{\varepsilon}||_V^2)du}||X_s^{n,\varepsilon}-X_s^{\varepsilon}||_V^2ds,
\end{eqnarray}
if $\gamma \geq 2C_1$, where $C_1$ is the constant appeared in (\ref{2.6}).
\vskip 0.3cm
Similar to the proofs of (\ref{5.02}), (\ref{5.03}), using Burkh\"older's inequality, we obtain from (\ref{5.34}), (\ref{5.34-1}) that for $t\leq T$,
\begin{eqnarray}\label{5.35}
&&E[\sup_{0\leq s\leq t}e^{-\gamma \int_0^s(||X_u^{n,\varepsilon}||_V^2+||X_u^{\varepsilon}||_V^2)du}|X_s^{n,\varepsilon}-X_s^{\varepsilon}|_H^2]\nonumber\\
&&+E[\int_0^te^{-\gamma \int_0^s(||X_u^{n,\varepsilon}||_V^2+||X_u^{\varepsilon}||_V^2)du}||X_s^{n,\varepsilon}-X_s^{\varepsilon}||_V^2ds]\nonumber\\
&\leq & \frac{1}{4}E[\sup_{0\leq s\leq t}e^{-\gamma \int_0^s(||X_u^{n,\varepsilon}||_V^2+||X_u^{\varepsilon}||_V^2)du}
|X_s^{n,\varepsilon}-X_s^{\varepsilon}|_H^2]\nonumber\\
&+& CE[ \int_0^te^{-\gamma \int_0^s(||X_u^{n,\varepsilon}||_V^2+||X_u^{\varepsilon}||_V^2)du}
|X_s^{n,\varepsilon}-X_s^{\varepsilon}|_H^2 ds]\nonumber\\
&+&CE[\int_0^t\int_{|x|\leq \varepsilon}e^{-\gamma \int_0^s(||X_u^{n,\varepsilon}||_V^2+||X_u^{\varepsilon}||_V^2)du}|\frac{1}{\alpha(\epsilon)}(\sigma_n(X_{s-}^{n,\varepsilon})x-\sigma_n(X_{s-}^{\varepsilon})x)|_H^2 ds \nu (dx)]\nonumber\\
&+&CE[\int_0^t\int_{|x|\leq \varepsilon}e^{-\gamma \int_0^s(||X_u^{n,\varepsilon}||_V^2+||X_u^{\varepsilon}||_V^2)du}|\frac{1}{\alpha(\epsilon)}(\sigma_n(X_{s-}^{\varepsilon})x-\sigma(X_{s-}^{\varepsilon})x)|_H^2 ds \nu (dx)]\nonumber\\
&\leq & \frac{1}{4}E[\sup_{0\leq s\leq t}e^{-\gamma \int_0^s(||X_u^{n,\varepsilon}||_V^2+||X_u^{\varepsilon}||_V^2)du}|X_s^{n,\varepsilon}-X_s^{\varepsilon}|_H^2]\nonumber\\
&+& CE[ \int_0^te^{-\gamma \int_0^s(||X_u^{n,\varepsilon}||_V^2+||X_u^{\varepsilon}||_V^2)du}|X_s^{n,\varepsilon}-X_s^{\varepsilon}|_H^2) ds]\nonumber\\
&+&CE[\int_0^t\int_{|x|\leq \varepsilon}e^{-\gamma \int_0^s(||X_u^{n,\varepsilon}||_V^2+||X_u^{\varepsilon}||_V^2)du}|\frac{1}{\alpha(\epsilon)}(\sigma_n(X_{s-}^{\varepsilon})x-\sigma(X_{s-}^{\varepsilon})x)|_H^2 ds \nu (dx)],\nonumber\\
&&
\end{eqnarray}
 where uniform Lipschitz constant of $\sigma_n$ has been used. Applying the Gronwall's inequality we obtain
 \begin{eqnarray}\label{5.36}
&&E[\sup_{0\leq s\leq t}e^{-\gamma \int_0^s(||X_u^{n,\varepsilon}||_V^2+||X_u^{\varepsilon}||_V^2)du}|X_s^{n,\varepsilon}-X_s^{\varepsilon}|_H^2]\nonumber\\
&&+E[\int_0^te^{-\gamma \int_0^s(||X_u^{n,\varepsilon}||_V^2+||X_u^{\varepsilon}||_V^2)du}||X_s^{n,\varepsilon}-X_s^{\varepsilon}||_V^2ds]\nonumber\\
&\leq &CE[\int_0^T|\sigma_n(X_{s}^{\varepsilon})-\sigma(X_{s}^{\varepsilon})|_H^2 ds].
\end{eqnarray}
For any $M>0$, we have
\begin{eqnarray}\label{5.37}
&&E[\int_0^T|\sigma_n(X_{s}^{\varepsilon})-\sigma(X_{s}^{\varepsilon})|_H^2 ds]\nonumber\\
&=&E[\int_0^T|\sigma_n(X_{s}^{\varepsilon})-\sigma(X_{s}^{\varepsilon})|_H^2 ds, \sup_{0\leq s\leq T}|X_s^{\varepsilon}|_H\leq M]\nonumber\\
&+&E[\int_0^T|\sigma_n(X_{s}^{\varepsilon})-\sigma(X_{s}^{\varepsilon})|_H^2 ds, \sup_{0\leq s\leq T}|X_s^{\varepsilon}|_H>M]\nonumber\\
&\leq& T\sup_{|z|\leq M}|\sigma_n(z)-\sigma(z)|_H^2+C T\frac{1}{M}(1+E[ \sup_{0\leq s\leq T}|X_s^{\varepsilon}|_H^3])\nonumber\\
&\leq& \sup_{|z|\leq M}|\sigma_n(z)-\sigma(z)|_H^2+C T\frac{1}{M},
\end{eqnarray}
where (\ref{5.0}) has been used. Since $M$ can be chosen as large as we wish, together with (\ref{5.36}) and (H.3)(ii) we deduce that
\begin{eqnarray}\label{5.37-1}
&&\lim_{n\rightarrow \infty}\sup_{\varepsilon}E[\sup_{0\leq s\leq T}e^{-\gamma \int_0^s(||X_u^{n,\varepsilon}||_V^2+||X_u^{\varepsilon}||_V^2)du}|X_s^{n,\varepsilon}-X_s^{\varepsilon}|_H^2]\nonumber\\
&=&0.
\end{eqnarray}
\vskip 0.3cm
For any given $\delta_1>0$, in view of (\ref{5.33-1}), (\ref{5.33-2}), we can choose a positive constant $M_1$ such that
\begin{eqnarray}\label{5.37-2}
&&\sup_{n,\varepsilon} P(\sup_{0\leq t\leq T}|X_t^{n,\varepsilon}-X_t^{\varepsilon}|_H>\delta, \int_0^T(||X_s^{n,\varepsilon}||_V^2+||X_s^{\varepsilon}||_V^2)ds>M_1)\nonumber\\
&\leq &\sup_{n,\varepsilon} P(\int_0^T(||X_s^{n,\varepsilon}||_V^2+||X_s^{\varepsilon}||_V^2)ds>M_1)\leq \frac{\delta_1}{2}.
\end{eqnarray}
On the other hand,
\begin{eqnarray}\label{5.37-3}
&&\sup_{\varepsilon} P(\sup_{0\leq t\leq T}|X_t^{n,\varepsilon}-X_t^{\varepsilon}|_H^2>\delta, \int_0^T(||X_s^{n,\varepsilon}||_V^2+||X_s^{\varepsilon}||_V^2)ds\leq M_1)\nonumber\\
&\leq &\sup_{\varepsilon} P(\sup_{0\leq s\leq T}e^{-\gamma \int_0^s(||X_u^{n,\varepsilon}||_V^2+||X_u^{\varepsilon}||_V^2)du}|X_s^{n,\varepsilon}-X_s^{\varepsilon}|_H^2\geq e^{-\gamma M_1}\delta^2)\nonumber\\
&\leq &e^{\gamma M_1}\frac{1}{\delta^2}\sup_{\varepsilon}E[\sup_{0\leq s\leq T}e^{-\gamma \int_0^s(||X_u^{n,\varepsilon}||_V^2+||X_u^{\varepsilon}||_V^2)du}|X_s^{n,\varepsilon}-X_s^{\varepsilon}|_H^2].
\end{eqnarray}
It follows from (\ref{5.37-1}) and (\ref{5.37-3}) that there exists $N>0$ such that for $n\geq N$,
\begin{eqnarray}\label{5.37-4}
&&\sup_{n,\varepsilon} P(\sup_{0\leq t\leq T}|X_t^{n,\varepsilon}-X_t^{\varepsilon}|_H>\delta, \int_0^T(||X_s^{n,\varepsilon}||_V^2+||X_s^{\varepsilon}||_V^2)ds\leq M_1)\nonumber\\
&\leq &\frac{\delta_1}{2}.
\end{eqnarray}
Combining (\ref{5.37-2}) and (\ref{5.37-4}) together yields (\ref{5.32})

Finally we prove that $\mu^{\varepsilon}$ converges to $\mu$. Let $\mu^{\varepsilon}_n$, $\mu_n$ denote respectively the laws of $X^{n,\varepsilon}$ and $X^{n}$ . Let $G$ be a bounded, uniformly continuous function on $E:=D([0,T], H)$. For any $n\geq 1$,  we write
\begin{eqnarray}\label{5.38}
&&\int_EG(w)\mu^{\varepsilon}(dw)-\int_EG(w)\mu(dw)\nonumber\\
&=&\int_EG(w)\mu^{\varepsilon}(dw)-\int_EG(w)\mu^{\varepsilon}_n(dw)+\int_EG(w)\mu_n^{\varepsilon}(dw)-\int_EG(w)\mu_n(dw)\nonumber\\
&& \quad\quad \quad +\int_EG(w)\mu_n(dw)-\int_EG(w)\mu(dw)\nonumber\\
&=& E[G(X^{\varepsilon})-G(X^{n,\varepsilon})]+( \int_EG(w)\mu_n^{\varepsilon}(dw)-\int_EG(w)\mu_n(dw))\nonumber\\
&&\quad\quad +E[G(X^{n})-G(X)].
\end{eqnarray}
Give any $\delta>0$. Since $G$ is uniformly continuous, there exists $\delta_1>0$ such that
\begin{equation}\label{5.40}
|E[(G(X^{\varepsilon})-G(X^{n,\varepsilon})), \sup_{0\leq s\leq T}|X_s^{n,\varepsilon}-X_s^{\varepsilon}|_H\leq \delta_1]|\leq \frac{\delta}{4}
\end{equation}
for all $n\geq 1, \varepsilon>0$. In view of (\ref{5.32}) and (\ref{5.33}), there exists $N_1$,
\begin{eqnarray}\label{5.41}
&&\sup_{\varepsilon}|E[(G(X^{\varepsilon})-G(X^{N_1,\varepsilon})), \sup_{0\leq s\leq T}|X_s^{N_1,\varepsilon}-X_s^{\varepsilon}|_H> \delta_1]|\nonumber\\
&\leq& C\sup_{\varepsilon}P(\sup_{0\leq s\leq T}|X_s^{N_1,\varepsilon}-X_s^{\varepsilon}|_H>\delta_1)  \leq \frac{\delta}{4},
\end{eqnarray}
and
\begin{equation}\label{5.42}
|E[(G(X^{N_1})-G(X))]|\leq \frac{\delta}{4}.
\end{equation}
On the other hand, by Theorem 3.4, there exists $\varepsilon_1>0$ such that for $\varepsilon \leq \varepsilon_1$,
\begin{equation}\label{5.42}
|\int_EG(w)\mu_{N_1}^{\varepsilon}(dw)-\int_EG(w)\mu_{N_1}(dw))|\leq \frac{\delta}{4}.
\end{equation}
Putting (\ref{5.38})---(\ref{5.42}) together we obtain that for $\varepsilon \leq \varepsilon_1$,
$$|\int_EG(w)\mu^{\varepsilon}(dw)-\int_EG(w)\mu(dw)|\leq \delta.$$
Since $\delta>0$ is arbitrarily small, we deduce that
$$\lim_{\varepsilon\rightarrow 0}\int_EG(w)\mu^{\varepsilon}(dw)=\int_EG(w)\mu(dw)$$
finishing the proof of the Theorem.$\blacksquare$

\begin{example}
Approximations of stochastic Burgers equations
\end{example}
Consider the stochastic Burgers equations on $[0,1]$:
\begin{eqnarray}\label{5.43}
du(t,\xi)=\frac{\partial^2}{\partial \xi^2}u(t,\xi) dt+\frac{1}{2}\frac{\partial}{\partial \xi}[u^2(t,\xi)]dt+\sigma(u(t,\xi))dB_t ,\\
u(t,0)=u(t,1)=0,t>0,
\end{eqnarray}
\begin{eqnarray}\label{5.44}
du^{\varepsilon}(t,\xi)=&\frac{\partial^2}{\partial \xi^2}u^{\varepsilon}(t,\xi) dt+\frac{1}{2}\frac{\partial}{\partial \xi}[(u^{\varepsilon})^2(t,\xi)]dt\\
&+\frac{1}{\alpha(\epsilon)}\int_{|x|\leq \varepsilon}\sigma(u^{\varepsilon}(t-, \xi))x\tilde{N}(dt,dx),\\
&u^{\varepsilon}(t,0)=u^{\varepsilon}(t,1)=0,t>0,
\end{eqnarray}
where $\sigma(\cdot)$ is a Lipschitz continuous function with $\sigma(0)=0$.

Let $V=H^1_0(0,1)$ with the norm
$$||v||_V:=\big(\int_0^1 (\frac{\partial u(\xi)}{\partial \xi})^2\,d\xi\big)^{\frac{1}{2}}=||v||.$$
Let $H:=L^2(0,1)$ be the $L^2$-space with inner product $(\cdot)$.\\
Set
$$Au=-\frac{\partial^2}{\partial \xi^2}u(\xi), \forall u \in H^2(0,1)\cap V.$$
Define for $k\geq 1$,
$$e_k(\xi)=\sqrt{2}sin(k\pi\xi), \xi\in [0,1].$$
Then $e_k, k\geq 1$ are eigenvectors of the operator $A$ with eigenvalues $\lambda_k=\pi^2k^2$, which forms an orthonormal basis of
the Hilbert space $H$. For $u\in V$, define
$$ B(u):= u(\xi)\frac{\partial}{\partial \xi}u(\xi), \quad \sigma(u):=\sigma(u(\xi)).$$
By the Lipschitz continuity of $\sigma$, it is easily seen that
\begin{equation}\label{5.45}
||\sigma(u)||_V\leq C(1+||u||_V).
\end{equation}
Now let us show that  $B(u)$ satisfies the condition (H.2). First it holds that
$$<B(u),u>=\int_0^1 u^2(\xi)\frac{\partial}{\partial \xi}u(\xi)d\xi=\frac{1}{3}[u^3(1)-u^3(0)]=0.$$
Note that $\bar{e}_k=\frac{1}{\sqrt{\lambda_k}}e_k, k\geq 1$ forms an orthonomal basis of $V$. Thus, for $u\in V$, we have
\begin{eqnarray}\label{5.46}
||B(u)||_{V^*}^2&=&\sum_{k=1}^{\infty}<B(u),\bar{e}_k>^2 \nonumber\\
&=&\sum_{k=1}^{\infty}\left( \frac{1}{2}\frac{1}{\sqrt{\lambda_k}}\int_0^1\frac{\partial}{\partial \xi}[u^2(\xi)]e_k(\xi)d\xi\right)^2\nonumber\\
&=&\sum_{k=1}^{\infty}\left( \frac{1}{2}\frac{1}{\sqrt{\lambda_k}}\int_0^1u^2(\xi)\frac{\partial}{\partial \xi}e_k(\xi)d\xi\right)^2\nonumber\\
&=&\sum_{k=1}^{\infty}\left (\frac{1}{2} \int_0^1u(\xi)^2\sqrt{2}cos(k\pi\xi)d\xi\right)^2 \nonumber\\
&\leq& C \int_0^1u^{\varepsilon}(s,\xi)^4d\xi=C|u|_{L^4}^4,
\end{eqnarray}
where we have used the fact that $\{ \sqrt{2}cos(k\pi\xi); k\geq 1\}$ also forms an orthonormal system of $L^2(0,1)$. Using the following well known
interpolation inequality
\begin{equation}\label{5.47}
|u|_{L^4}^4\leq C |u|_{H}^{3}||u||_{V},
\end{equation}
we obtain from (\ref{5.46}) that
$$||B(u)||_{V^*}\leq C |u|_{H}^{\frac{3}{2}}||u||_{V}^{\frac{1}{2}}$$
proving (H.2)(iii) with $\gamma=\frac{1}{2}$. Finally we will check (H.2)(ii). Let $u,v\in V$. We have
\begin{eqnarray}\label{5.47}
&&<B(u)-B(v), u-v>=\frac{1}{2}\int_0^1\frac{\partial}{\partial \xi}[u^2(\xi)-v^2(\xi)](u(\xi)-v(\xi))d\xi\nonumber\\
&=&-\frac{1}{2}\int_0^1(u^2(\xi)-v^2(\xi))\frac{\partial}{\partial \xi}(u(\xi)-v(\xi))d\xi\nonumber\\
&\leq &\frac{1}{2}\int_0^1(\frac{\partial}{\partial \xi}(u(\xi)-v(\xi)))^2d\xi\nonumber\\
&&+C\int_0^1(u(\xi)-v(\xi))^2(u(\xi)+v(\xi))^2d\xi\nonumber\\
&\leq &\frac{1}{2}||u-v||_V^2+C|u-v|_H^2(||u||_{\infty}^2+||v||_{\infty}^2)\nonumber\\
&\leq& \frac{1}{2}||u-v||_V^2+C|u-v|_H^2(||u||_{V}^2+||V||_{V}^2),
\end{eqnarray}
which is (H.2)(ii).
\vskip 0.3cm
Now we can apply Theorem 3.5 to obtain the following convergence of the solutions of stochastic Burgers equations.
\begin{theorem}
Let $u^{\varepsilon}, u$ be solutions to the stochastic Burgers equations (\ref{5.44}) and (\ref{5.43}). Then $u^{\varepsilon}$
converges weakly to $u$ in the space $D([0,T]; H)$.
\end{theorem}

\end{document}